\documentclass[11pt]{article}
\usepackage{dsfont}
\usepackage{graphics}
\usepackage[demo]{graphicx}
\usepackage{epsfig}
\usepackage{pstricks}
\usepackage[normal]{subfigure}
\usepackage[latin1]{inputenc}
\usepackage[english,french]{babel}
\usepackage{relsize,exscale}
\usepackage{makeidx}
\usepackage{enumitem}
\usepackage{amsfonts,amssymb,amsmath}
\usepackage{graphicx}
\usepackage{color}
\usepackage{multirow}
\usepackage{mathrsfs}
\usepackage[normalem]{ulem}
\usepackage{cancel}
\usepackage{bbm}
\usepackage{color}
\definecolor{dred}{rgb}{0.92,0,0}
\definecolor{dgreen}{rgb}{0,0.92,0}
\definecolor{dblue}{rgb}{0,0,0.92}
\definecolor{dyellow}{rgb}{0.95,0.95,0}

\usepackage{multirow}
\newcommand{\hs}{\hspace{0.1cm}}

\newcommand{\sa}{\\ [0.2cm]}
\newcommand{\R}{\mathbb{R}}
\newcommand{\N}{\mathbb{N}}

\newcommand{\be}{\begin{equation}}
\newcommand{\ee}{\end{equation}}
\newcommand{\bes}{\begin{equation*}}
\newcommand{\ees}{\end{equation*}}
\newcommand{\bea}{\begin{eqnarray}}
\newcommand{\eea}{\end{eqnarray}}
\newcommand{\beas}{\begin{eqnarray*}}
\newcommand{\eeas}{\end{eqnarray*}}

\makeatletter
\newcommand*{\rom}[1]{\expandafter\@slowromancap\romannumeral #1@}

\def\D{\displaystyle}
\def\ds{\displaystyle}

\newtheorem{theorem}{Theorem}[section]
\newtheorem{lemma}[theorem]{Lemma}

\newtheorem{e-definition}[theorem]{Definition\rm}

\graphicspath{
{../test-numeriques-FF++/function-de-Runge/}
{./Figures/}
{../Figures/}
{./}
}

\title{Numerical validation of probabilistic laws\\
 to evaluate finite element error estimates
 }
\author{Jo\"el Chaskalovic \thanks{D'Alembert,
Sorbonne University, Paris, France, (\emph{email}: jch1826@gmail.com)}
\qquad
Franck Assous
\thanks{
Department of Mathematics, Ariel University, 40700 Ariel, Isra\"el, (\emph{email}: franckassous55@gmail.com).}
}
\date{}

\begin{document}

\maketitle
\selectlanguage{english}
\begin{abstract}
\noindent We propose a numerical validation of a probabilistic approach applied to estimate the relative accuracy between two Lagrange finite elements $P_k$ and $P_m, (k<m)$. In particular, we show practical cases where finite element $P_{k}$ gives more accurate results than finite element $P_{m}$. This illustrates the theoretical probabilistic framework we recently derived in order to evaluate the actual accuracy. This also highlights the importance of the extra caution required when comparing two numerical methods, since the classical results of error estimates concerns only the asymptotic convergence rate.

\end{abstract}
\noindent {\em keywords}: Numerical validation, Error estimates, Finite elements, Bramble-Hilbert lemma, Probability.
%
%
%
\section{Introduction}\label{intro}

\noindent Finite element methods and among them, error estimates play a significant role in the development of numerical  methods. Very often, the success of a numerical method depends on its performance in terms of efficiency and accuracy. For this reason, it is still an active subject of research, as observed, for instance, with the considerable interest received by the discontinuous Galerkin methods in the past decades; see e.g. an introduction for elliptic problems in \cite{ABCM02}, the book \cite{HeWa08} or the pioneering work \cite{LeRa74}.\sa
\noindent Since the seminal papers of Strang and Fix \cite{StFi73}, Ciarlet and Raviart  \cite{Ciarlet_Raviart},  Babuska \cite{Babu71} and Bramble and Hilbert \cite{BrHi70}, along with co-workers, a large amount of work has been published, the purpose of which was to derive and expand error estimates in different configurations. Here, we are concerned with {\em a priori} error estimates, that aim to find upper bounds for the error between the exact solution $u$ and its finite element approximation $u_h$. More precisely, these estimates describe how the finite element error $\|u-u_h\|$, for a given norm, goes to 0 with mesh size $h$ (i.e. the largest diameter of the elements in a given mesh). In addition, these estimates involve a constant, generally unknown, which leads to only get an upper bound for the approximation error.\sa
\noindent In addition, quantitative uncertainties do exist in finite element methods; these are based on the way the mesh grid generator creates the mesh which is used to compute the finite element approximation $u_h$, or since the equations are not exactly solved due to round-off errors. In previous papers \cite{CMAM_2019}, \cite{CMAM_2020}, we investigated the error resulting from a partial non-control of the mesh size. For this purpose, we have considered the approximation error as a random variable, and we have evaluated  the relative accuracy between two Lagrange finite elements with the help of a probabilistic approach. In the same way, one can find in \cite{Nova1}, \cite{Nova2} a probabilistic approach to evaluate error bounds in numerical analysis.\sa
\noindent In this work, we {\em numerically} study the {\em a priori error} estimate due to the discretization of a linear variational problem by a finite element method, using standard polynomials. Our aim is to compare the probabilistic laws we derived with statistical results, when two different degrees of the polynomials are used, for a fixed value of the mesh size. Since the effective dependence of the accuracy on the mesh size is a central question, it could help one to understand the saturation assumption that is often used in {\it a posteriori} error analysis \cite{OdAi00}. Indeed, we use here a probabilistic approach which differs from the methods involved in {\it a posteriori} error analysis. Nonetheless, we will show examples where $P_k$ finite element is more likely accurate than $P_m$, $k<m$, which can be related to the invalidity of the saturation assumption \cite{DoNo}.\sa
\noindent The paper is structured as follows: Section \ref{models} summarizes the results of \cite{CMAM_2019}, \cite{CMAM_2020}, which are necessary for one to understand the numerical experiments and their analysis. The main results are the geometrical interpretation of error estimates and the two probabilistic laws we deduced for finite element accuracy. Section \ref{NumericalResults} is devoted to the numerical results, which illustrate the new probabilistic way we propose to evaluate the accuracy between two finite elements. We basically consider two numerical problems, a stiff one and a smooth one, and we compare, for each of them, the behavior of the theoretical probabilistic models with the statistical results. Concluding remarks follow.
\section{Probabilistic models and finite elements accuracy}\label{models}
\subsection{Error estimates revisited}\label{sub1models}
\noindent Consider $\Omega$ an open bounded and non-empty subset of $\R^{n}$, and let $\Gamma$ denote its boundary, assumed to be $C^1- $piecewise. We also introduce an Hilbert space $V$ endowed with a norm, $ \left\|.\right\|_{V}$, and a bilinear, continuous and $V-$elliptic form $a(\cdot,\cdot)$  defined on $V \times V$. Finally, $l(\cdot)$ denotes a linear continuous form defined on~$V$.\sa
Let $u \in V$ be the unique solution to the second order elliptic variational formulation
\begin{equation}\label{VP}
\left\{
\begin{array}{l}
\mbox{Find } u \in   V \mbox{ solution to:} \\[0.1cm]
a(u,v) = l(v), \quad\forall v \in V\,.
\end{array}
\right.
\end{equation}
In this paper, we will focus on the simple case where $V$ is the usual Sobolev space of distributions $H^1(\Omega)$.  More general cases can be found in \cite{ChAs20}.\sa
Let us now introduce the finite-dimensional subspace $V_h$ of $V$, and consider $u_{h}\in V_{h}$ an approximation of $u$, solution to the approximate variational formulation
$$
\left\{
\begin{array}{l}
\mbox{Find } u_{h} \in   V_h \mbox{ solution to:} \\[0.1cm]
a(u_{h},v_{h}) = l(v_{h}),\quad \forall v_{h} \in V_h.
\end{array}
\right.
$$
In what follows, we are interested in evaluating error bounds for finite element methods. Hence, we first assume that domain $\Omega$ is exactly covered by a mesh ${\mathcal T}_h$ composed by $N_{s}$ n-simplexes $K_{j}, (1 \leq j \leq N_{s}),$ which respects classical rules of regular discretization, (see for example \cite{ChaskaPDE} for the bidimensional case, or \cite{RaTho82} in $\R^n$). We also denote by $P_k(K_{j})$ the space of polynomial functions defined on a given n-simplex $K_{j}$ of degree less than or equal to $k$, ($k \geq$ 1). \sa
Our study relies on the results of \cite{RaTho82}. Let $\|.\|_{1}$ be the classical norm in $H^1(\Omega)$ and $|.|_{k+1}$ the semi-norm in $H^{k+1}(\Omega)$, and let $h$ be the mesh size, namely the largest diameter of the elements of the mesh ${\mathcal T}_h$. We thus have:
\begin{lemma}\label{Thm_error_estimate}
Suppose that there exists an integer $k \geq 1$ such that the approximation $u_h$ of $V_h$ is a continuous piecewise function composed by polynomials which belong to $P_k(K_{j}), (1\leq j\leq  N_{s})$. \sa
Then, if the exact solution $u$ belongs to $H^{k+1}(\Omega)$, we have the following error estimate:
\begin{equation}\label{estimation_error}
\|u_h-u\|_{1} \hs \leq \hs \mathscr{C}_k\,h^k \, |u|_{k+1}\,,
\end{equation}
where $\mathscr{C}_k$ is a positive constant independent of $h$.
\end{lemma}

\noindent Now, let us consider two families of Lagrange finite elements $P_k$ and $P_m$ for two values $(k,m)\in \N^{*2}$, $(k < m)$. Assuming that the solution $u$ to (\ref{VP}) belongs to $H^{m+1}(\Omega)$,  inequality (\ref{estimation_error}) can be written as
\begin{eqnarray}
\|u^{(k)}_h-u\|_{1} \hs & \leq & \hs \mathscr{C}_k h^{k}\, |u|_{k+1}, \label{Constante_01} \\
\|u^{(m)}_h\hspace{-0.09cm}-u\|_{1} \hs & \leq & \hs \mathscr{C}_m h^{m}\, |u|_{m+1}\,, \label{Constante_02}
\end{eqnarray}
where $u^{(k)}_h$ and $u^{(m)}_h$ respectively denote the $P_k$ and $P_m$ Lagrange finite element approximations of $u$.\sa
In this article, following a series of previous papers \cite{CMAM_2019}-\cite{ChAs20} where a theoretical analysis was performed, we are interested in numerical applications. To this end, for a given mesh size $h$, two independent meshes for $P_k$ and $P_m$ are built by a mesh generator. Usually, one considers inequalities (\ref{Constante_01}) and (\ref{Constante_02}) so as to conclude that, when $h$ goes to zero, $P_m$ is more accurate that $P_k$, since $h^m$ goes faster to zero than $h^k$. \sa
However, in practical numerical applications, the size of the mesh is chosen according to the desired accuracy, so that $h$ has a fixed value. Consequently, this way of comparison is no more relevant. For this reason, we mean to identify the relative accuracy between $P_k$ and $P_m, (k<m)$, for a given value of $h$.
\subsection{Two probabilistic laws}\label{sub2models}
\begin{figure}[h]
\centering
\includegraphics[width=10.cm]{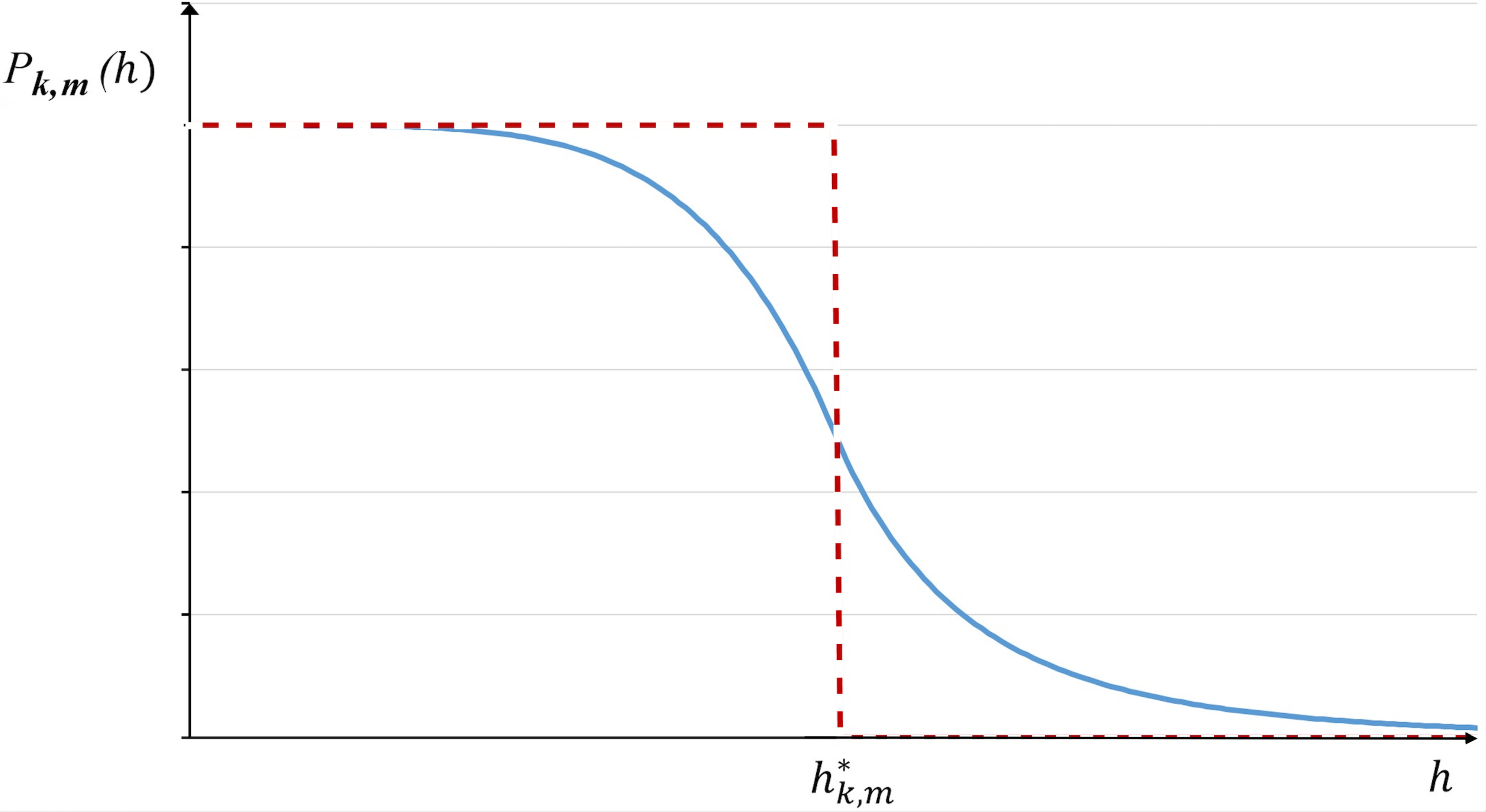}
\caption{Case $m-k\neq 1$: shape of the sigmoid distribution (\ref{Nonlinear_Prob}) (full line) and the two steps corresponding one (\ref{Heaviside_Prob}) (dashed line), $(P_{k,m}(h)\equiv Prob\{X^{(m)}(h) \leq X^{(k)}(h)\})$.} \label{Sigmoid}
\end{figure}

\noindent In \cite{CMAM_2019}-\cite{CMAM_2020},  we introduced a probabilistic approach that provides a coherent framework for modeling uncertainties in finite element approximations: such uncertainties may come from the way the meshes are created by computer algorithms, leading to a partial non-control of the mesh, even for a given maximum mesh size.\sa
In this framework, values $\|u^{(k)}_h-u\|_{1}$ and $\|u^{(m)}_h-u\|_{1}$ are viewed as two random variables, respectively denoted as $X^{(k)}(h)$ and $X^{(m)}(h)$, whose support is $\big[0,  \mathscr{C}_i |u|_{i+1} h^i\big], (i=k \mbox{ or } i=m),$ according to inequalities (\ref{Constante_01}) and (\ref{Constante_02}). Our goal is thus to derive a probabilistic law for the event
$$
\D\left\{X^{(m)}(h) \leq X^{(k)}(h)\right\} \equiv \left\{\|u^{(m)}_h-u\|_{1} \leq \|u^{(k)}_h-u\|_{1}\right\},
$$
which corresponds to the relative accuracy between finite elements $P_k$ and $P_m$. For this purpose, we first introduce the random events $A$ and $B$ defined by:
\begin{eqnarray*}
A  \equiv \left\{X^{(m)}(h) \leq X^{(k)}(h)\right\}, \, &&
B  \equiv  \left\{ X^{(k)}(h) \in \big[\mathscr{C}_m |u|_{m+1} h^m,\mathscr{C}_k |u|_{k+1} h^k\big]\right\}.
\end{eqnarray*}
Moreover, we proved in \cite{CMAM_2019} the following result:
\begin{lemma}\label{Two_Steps}
Let us assume that $A$ and $B$ are two independent events. Then, the probability law $P(A)$ of event $A$ is given by:
\vspace{-0.2cm}
\begin{equation}\label{Heaviside_Prob}
\D P(A)= \left |
\begin{array}{ll}
\hs 1 & \mbox{ if } \hs 0 < h < h^{*}_{k,m}, \medskip \\
\hs 0 & \mbox{ if } \hs h>  h^{*}_{k,m},
\end{array}
\right.
\end{equation}
where $h^{*}_{k,m}$ is defined by:
\begin{equation}\label{h*}
\D h^{*}_{k,m} \equiv\left( \frac{\mathscr{C}_k |u|_{k+1}}{\mathscr{C}_m |u|_{m+1}}\right)^{\frac{1}{m-k}}.
\end{equation}
\end{lemma}
The shape of the probabilistic distribution, called the two-steps model, is depicted in Fig.\ref{Sigmoid}. Basically, it expresses the fact that, for $h < h^{*}_{k,m}$, finite element $P_m$ is \emph{almost surely} more accurate than $P_k$, whereas for $h > h^{*}_{k,m}$, $P_k$ becomes \emph{almost surely} more accurate than $P_m$.\sa
To relax the independence assumption of events $A$ and $B$, we also derived a second probabilistic law based on the uniform distribution of the random variable $X^{(k)}(h)$ over $\big[0,\mathscr{C}_k |u|_{k+1} \,h^{k}\big]$. In this context, we proved in \cite{CMAM_2019} the following theorem:
\begin{theorem}\label{The_nonlinear_law}
Let us assume that $X^{(i)}(h), (i=k,m),$ are independent and uniformly distributed on $[0, \mathscr{C}_i |u|_{i+1} h^i]$. Then, the probability $P(A)$ of event $A$ is given by:
\begin{equation}\label{Nonlinear_Prob}
\D P(A)= \left |
\begin{array}{ll}
\D \hs 1 - \frac{1}{2}\!\left(\!\frac{\!\!h}{h^{*}_{k,m}}\!\right)^{\!\!m-k} & \mbox{ if } \hs 0 < h \leq h^{*}_{k,m}, \\[0.5cm]
\D \hs \frac{1}{2}\!\left(\!\frac{h^{*}_{k,m}}{\!\!h}\!\right)^{\!\!m-k} & \mbox{ if } \hs h \geq h^{*}_{k,m}.
\end{array}
\right.
\end{equation}
\end{theorem}
The shape of this law, called the sigmoid model, is also plotted in Fig.\ref{Sigmoid}. As one can see,  for $h>h^{*}_{k,m}$, $P(A) \leq 0.5$: in that case, finite element $P_{m}$ is \emph{probably} overqualified.\sa
The purpose of the next section is to propose numerical examples that illustrate and validate this probabilistic approach by comparing statistical frequencies and the corresponding probabilities determined by (\ref{Heaviside_Prob}) or (\ref{Nonlinear_Prob}).
\section{Numerical results}
\label{NumericalResults}
\noindent In this section, we will illustrate our probabilistic approach on numerical examples, by evaluating the relative accuracy of two Lagrange finite elements. We have intentionally chosen a simple, standard example, in order to help us  numerically check the relevance of the proposed probabilistic distributions.\sa
\noindent Hence, we consider the following classical elliptic problem, with obvious notations\sa
{\em Find $u \in H^1(\Omega)$ solution to}
\begin{equation}\label{edpnum}
\left\{
\begin{array}{c}
-\Delta u = q \mbox{ in } \Omega\,, \\
u =h \mbox{ on } \partial \Omega\,,
\end{array}
\right.
\end{equation}
where, for simplicity, domain $\Omega$ is the open unit square in $\mathbb{R}^2$: $\Omega=]0,1[\times]0,1[$. The associated variational formulation, which is analogous to (\ref{VP}), can be readily derived. According to the choice of $q$ and $h$, we will consider as examples a {\em stiff} problem, where the solution exhibits rapid variations or, alternatively, a {\em very smooth} problem.\sa
One of the main ingredients of the method is the computation of $h^{*}_{k,m}$, as defined by (\ref{h*}). As one will see, it will be evaluated using a {\em maximum likelihood estimator}; see for instance \cite{Rossi}.\sa
In our case, this principle is applied as follows: for a given finite element $P_{k}$, we consider a number $N$ of different meshes with the same (maximum) mesh size $h$. Then, we compute:
\begin{equation}
\label{estimator-Ck}
\max_{{N,h}}\ds\frac{\|u^{(k)}_h-u\|_{1}}{h^{k}}\,
\end{equation}
which constitutes, using estimate (\ref{Constante_01}), the maximum likelihood estimator for $\mathscr{C}_{k}|u|_{k+1}$. Indeed, due to inequality (\ref{Constante_01}), quantity $\frac{X^{k}(h)}{h^k}$ is also a uniform random variable whose support is $[0, \mathscr{C}_{k}|u|_{k+1}]$. \sa
Then, one can show \cite{LeJeune} that for a given uniform random variable $Y$ whose support is $[0, \theta]$, $\theta$ being an unknown real parameter, the maximum likelihood estimator $\hat{\theta}$ is given by:
$$
\hat{\theta} = \max(Y_1,\dots,Y_N),
$$
where $(Y_1,\dots,Y_N)$ is a sample built with independent and identically distributed random variables $(Y_i)_{i=1,N}$, with the same distribution as $Y$.\sa
In our case, this implies that (\ref{estimator-Ck}) is the maximum likelihood estimator for $\mathscr{C}{k}|u|_{k+1}$, since $N$ and $h$ each take a finite number of values.\sa
Doing the same for another finite element $P_{m}$, we obtain that the estimator for $h^{*}_{k,m}$, denoted $\widehat{\,h^{*}_{k,m}}$, is defined by:
\begin{equation}\label{hstar_estime}
\widehat{\,h^{*}_{k,m}}=\left(\ds\frac{\ds\max_{{N,h}}\ds\frac{\|u^{(k)}_h-u\|_{1}}{h^{k}}}{\ds\max_{{N,h}}\ds\frac{\|u^{(m)}_h-u\|_{1}}{h^{m}}}\right)^{1/m-k}
\end{equation}
Then, one can easily compute the two probability laws introduced in subsection \ref{sub2models}. Indeed, as soon as $\widehat{\,h^{*}_{k,m}}$ is computed, functions (\ref{Heaviside_Prob}) and (\ref{Nonlinear_Prob}) are operational by replacing $h^{*}_{k,m}$ by $\widehat{\,h^{*}_{k,m}}$. All the numerical results below are computed in this way.\sa
In order to numerically check the validity of each model, we now compare the two probabilistic laws defined by (\ref{Heaviside_Prob}) and (\ref{Nonlinear_Prob}) with the corresponding statistical frequencies computed on the $N$ meshes, for each fixed value $h$ of the mesh size.\sa
To that end, we consider for two finite elements $P_{k}$ and $P_{m}$ ($k<m$), the same number $N$ of different meshes with the same (maximum) mesh size $h$. From there,  we compute the approximate solution $u^{(m)}_h$ and $u^{(k)}_h$, and we test if $\|u^{(m)}_h-u\|_{1} \leq \|u^{(k)}_h-u\|_{1}$. Then, we repeat the same process for different values of $h$, either lower or greater than  $\widehat{\,h^{*}_{k,m}}$. This gives, as a function of $h$, the percentage of cases where the approximation error of $P_m$ is lower than the approximation error of $P_k$. In all cases, we use package FreeFem++~\cite{Hech12} to compute the $P_k$ and $P_m$ finite element approximations.\sa
In the next subsection, we consider a stiff case, whereas in the following one, we deal with a very smooth example.
\subsection{A first stiff case}
To introduce such a stiff case, we consider the well-known Runge function $\varphi(t)=\ds\frac{1}{1+\alpha t^2}$ which takes $\alpha$ as a parameter, the classical Runge function corresponding to $\alpha=25$ (see \cite{Rung01}, \cite{Eppe87}). \sa
Since we first aim at building an exact solution $u(x,y)$ for (\ref{edpnum}), we consider solutions of the form $u(x,y)=f(x) g(y)$, where both $f(x)$ and $g(y)$ are Runge functions of parameter $\alpha$.\sa
To compute the derivatives of $u(x,y)$, we basically need the derivatives of the Runge function $f$. After some elementary algebra, we obtain the derivatives of $f(t)$ (namely $f'(t)=\ds\frac{-2 \alpha t}{(1+\alpha t^2)^2}, f''(t)=\ds\frac{2 \alpha (3 \alpha t^2 -1)}{(1+\alpha t^2)^3}$), from which the Laplacian of $u(x,y)$  can easily be derived. Indeed, by computing the second order partial derivatives $u_{xx}, u_{yy}$, we find that
\begin{equation}
\label{Rhs}
-\Delta u  = -(f''(x) g(y) +f(x) g''(y)) =
-\frac{2 \alpha (3 \alpha x^2 -1)}{(1+\alpha x^2)^3} \ds\frac{1}{1+\alpha y^2}-\ds\frac{1}{1+\alpha x^2}\frac{2 \alpha (3 \alpha y^2 -1)}{(1+\alpha y^2)^3}\,.
\end{equation}
We now set the right-hand side $q(x,y)$ of (\ref{edpnum}) equal to expression (\ref{Rhs}) above, so that
\begin{equation}
\label{uexact}
u(x,y)=\ds\left(\frac{1}{1+\alpha x^2}\right)\ds\left(\frac{1}{1+\alpha y^2}\right)
\end{equation}
is the exact solution of (\ref{edpnum}), provided that the Dirichlet boundary condition $h(x,y)$ is taken as the trace of $u(x,y)$ on the boundary $\partial \Omega$, that is
$$
\left\{
\begin{array}{ll}
h(x,0) =\ds\frac{1}{1+\alpha x^2}\,,  & h(0,y) = \ds\frac{1}{1+\alpha y^2}, \\
h(x,1) = \ds\frac{1}{1+\alpha x^2}\ds\frac{1}{1+\alpha} \,, & h(1,y) = \ds\frac{1}{1+\alpha}\ds\frac{1}{1+\alpha y^2}.
\end{array}
\right.
$$
In what follows, we analyze the relative accuracy between two Lagrange finite elements, $u(x,y)$, as defined in (\ref{uexact}), being the reference solution for comparison.
\subsubsection{$P_2$-$P_3$ comparison and $\alpha$-independence}
\begin{figure}[htb]
  \centering
  \includegraphics[width=8.5cm]{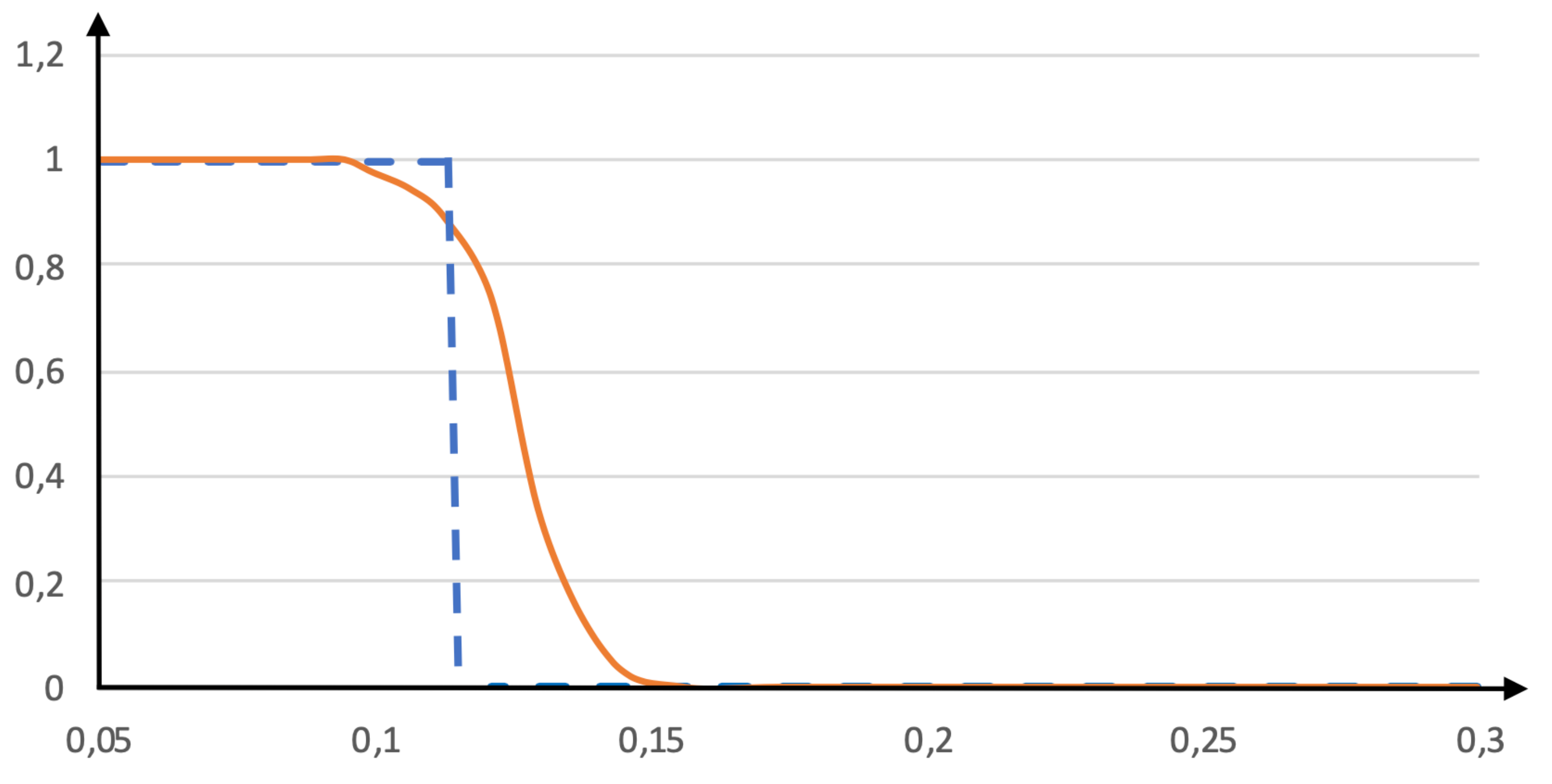}
  \caption{$P_2$ versus $P_3$ for the Runge function with $\alpha=500$. Comparison between the statistical frequencies (full line) and the probabilistic law (\ref{Heaviside_Prob}) (dotted line). 500 meshes are used for each value of $h$.} \label{P2P2alpha-500}
\end{figure}

\begin{figure}[htbp!]
\begin{tabular}{lr}
{
\hspace*{-0.5cm}
\includegraphics[width=8.5cm]{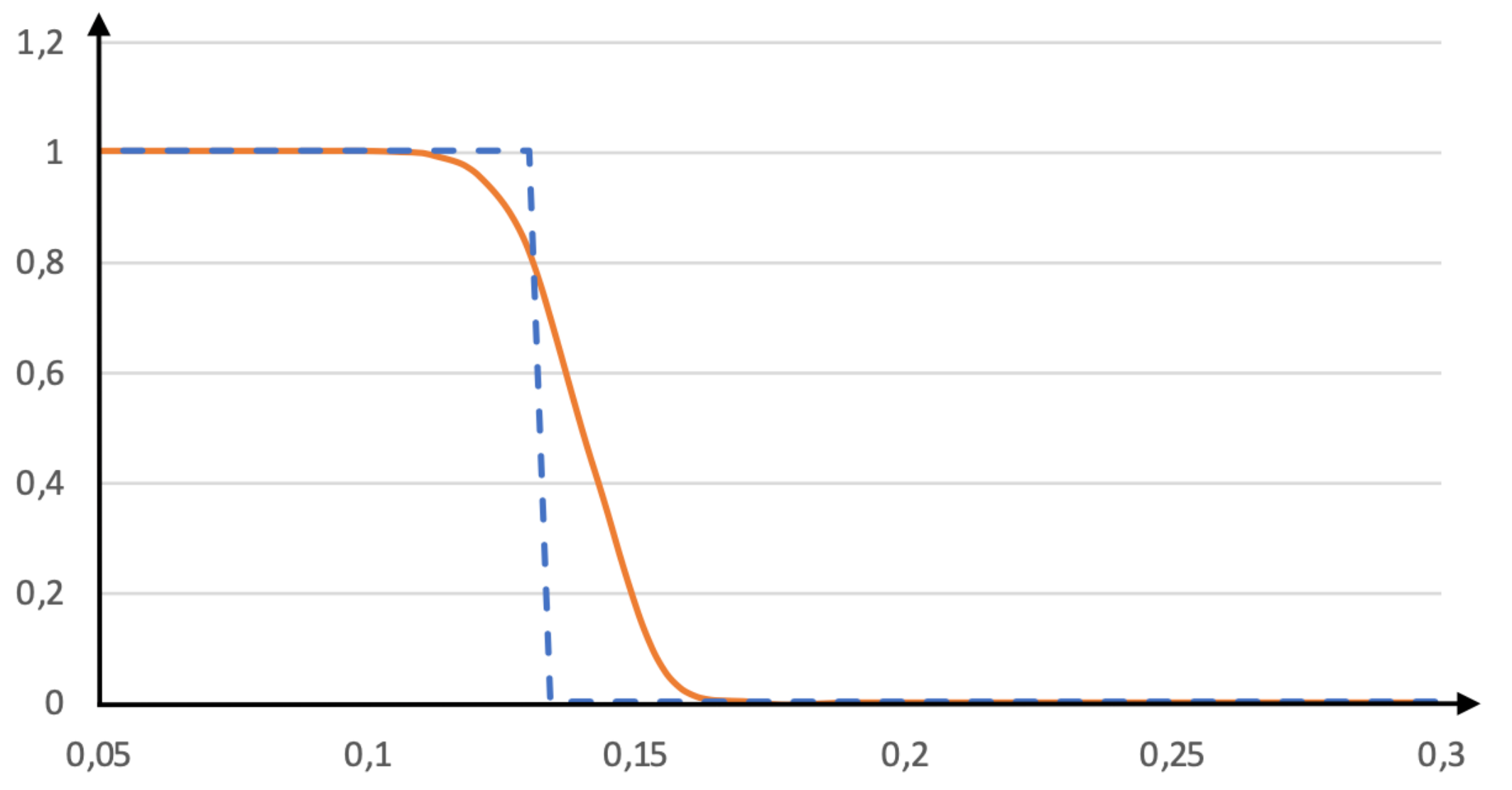}
}
&
\hspace*{-0.5cm}
{
\includegraphics[width=8.5cm]{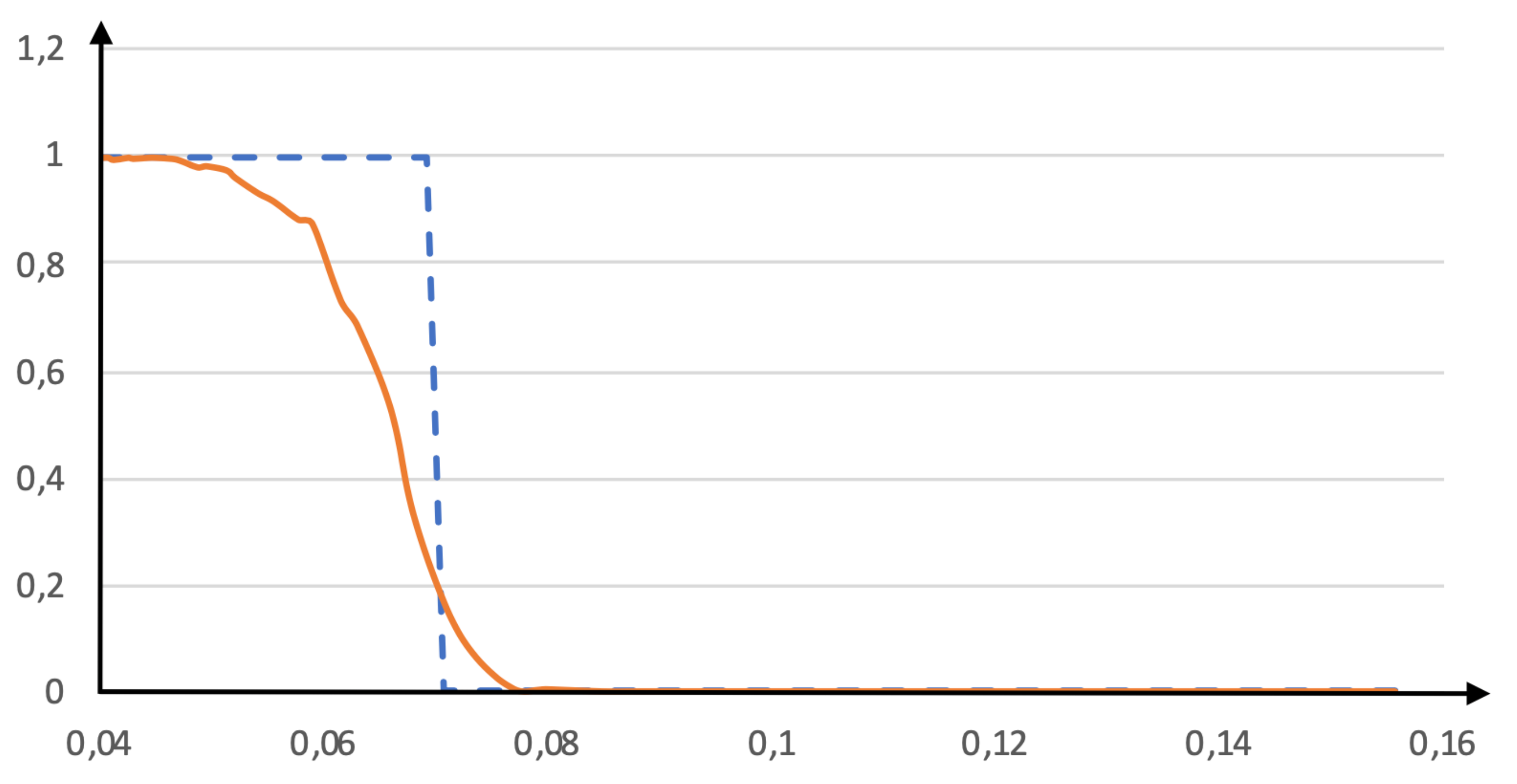}
}
\end{tabular}
\caption{$P_2$ versus $P_3$ for the Runge function with $\alpha=25$ (left) and $\alpha=2000$ (right). Statistical frequencies (full) and probabilistic law (\ref{Heaviside_Prob}) (dotted). 500 meshes are used for each value of $h$.} \label{P2P2alpha-25-2000}
\end{figure}

\noindent The first numerical test we present is devoted to a comparison between finite elements $P_2$ and $P_3$. We first choose $\alpha=500$. In that case, as explained above, we computed value $\widehat{\,h^{*}_{2,3}}$ as defined in (\ref{hstar_estime}) and obtained $\widehat{\,h^{*}_{2,3}} \simeq 0.12$.\sa
For this example, we have used values of $h$ varying from $0.05$ to $0.18$, and for each $h$ we have constructed $N=500$ different meshes with the same value of $h$. In Fig. \ref{P2P2alpha-500} we plot, on the same picture, the results obtained for the statistical frequencies (full line) and for the two-steps probability law (\ref{Heaviside_Prob}) (dotted line), as a function of $h$.
\noindent We then checked that the results do not depend on the value of $\alpha$. For this purpose, we repeated the same numerical experiments for $\alpha=25$ and $\alpha=2000$. The results, depicted in Fig. \ref{P2P2alpha-25-2000}, show the same behavior as previously. Of course, the value of $\widehat{\,h^{*}_{k,m}}$ does depend on $\alpha$ and we have $\widehat{\,h^{*}_{2,3}}\simeq 0,13$ for $\alpha=25$ and $\widehat{\,h^{*}_{2,3}}\simeq 0.07$ for $\alpha=2000$.\sa
\subsubsection{Comparison with $P_4$ finite element}
\begin{figure}[htbp!]
\begin{tabular}{lr}
{
\hspace*{-0.5cm}
\includegraphics[width=8.cm]{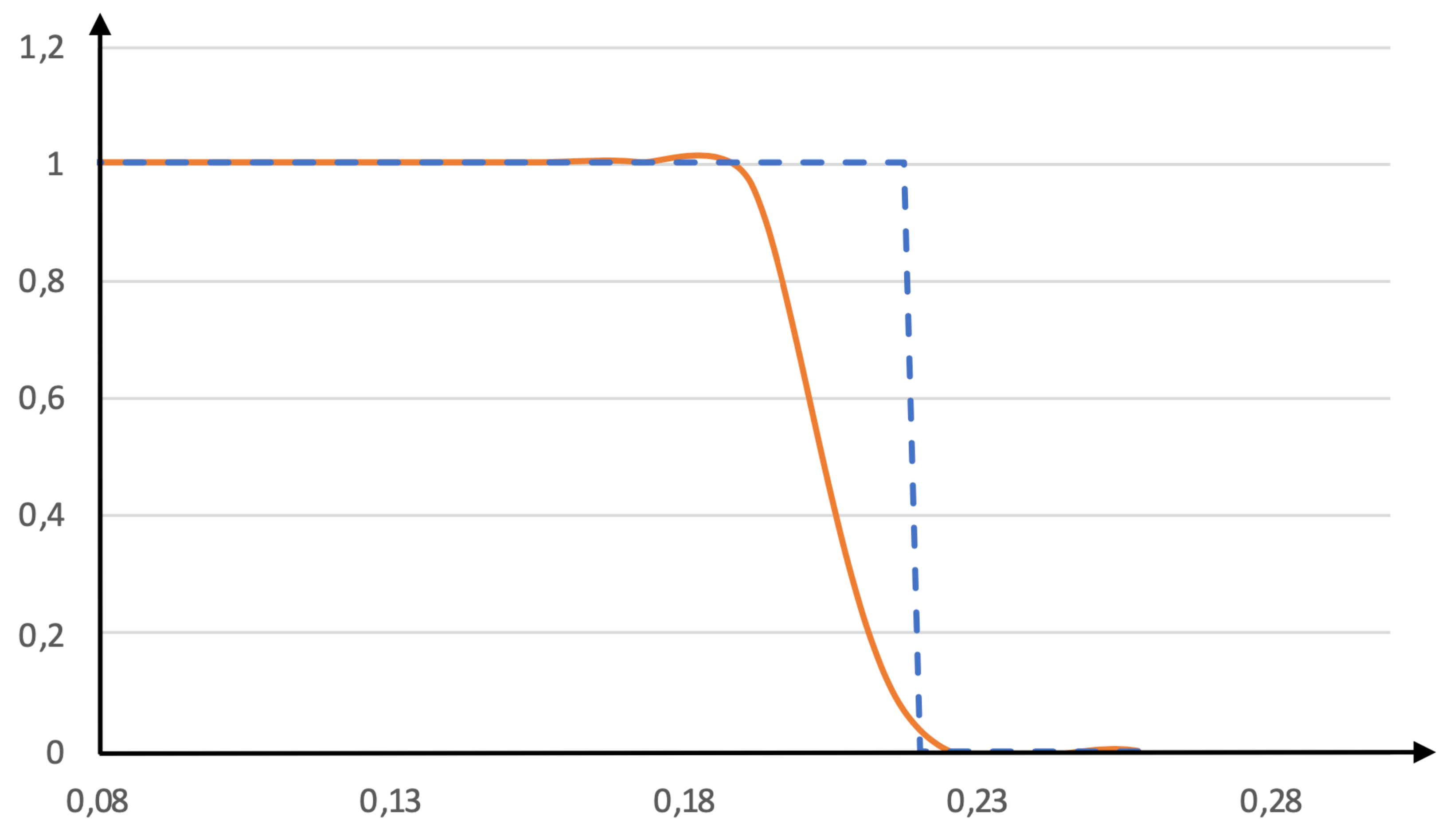}
}
&
\hspace*{-0.5cm}
{
\includegraphics[width=8.cm]{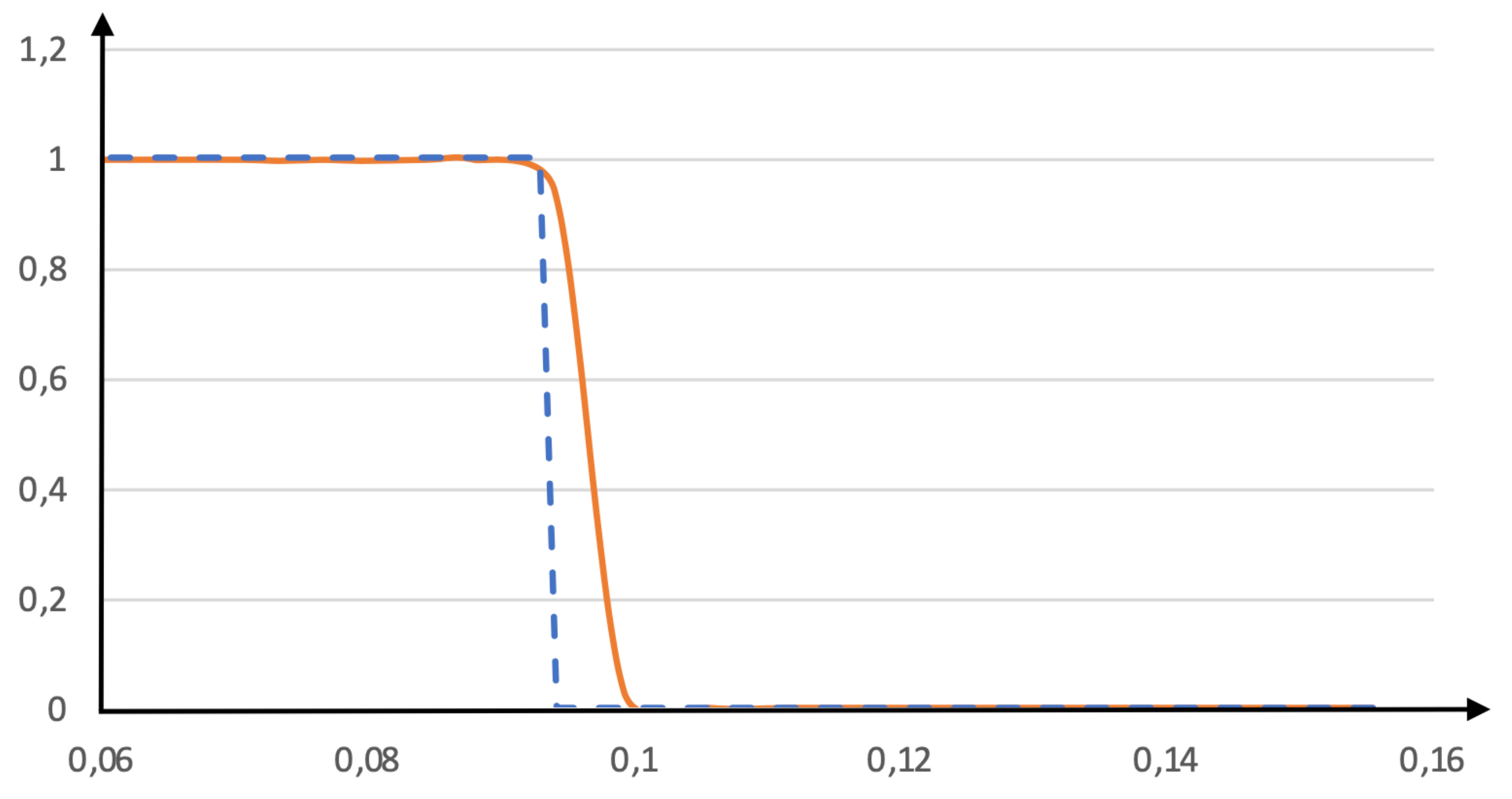}
}
\end{tabular}
\caption{Comparisons $P_3$-$P_4$ (left) and $P_2$-$P_4$ (right) for the Runge function with $\alpha=2000$. Statistical frequencies (full) and probabilistic law (\ref{Heaviside_Prob}) (dotted). 500 meshes are used for each value of $h$.} \label{P3P4-P2P4-alpha-2000}
\end{figure}

\begin{figure}[htbp!]
\begin{tabular}{lr}
{
\hspace*{-0.5cm}
\includegraphics[width=9.cm]{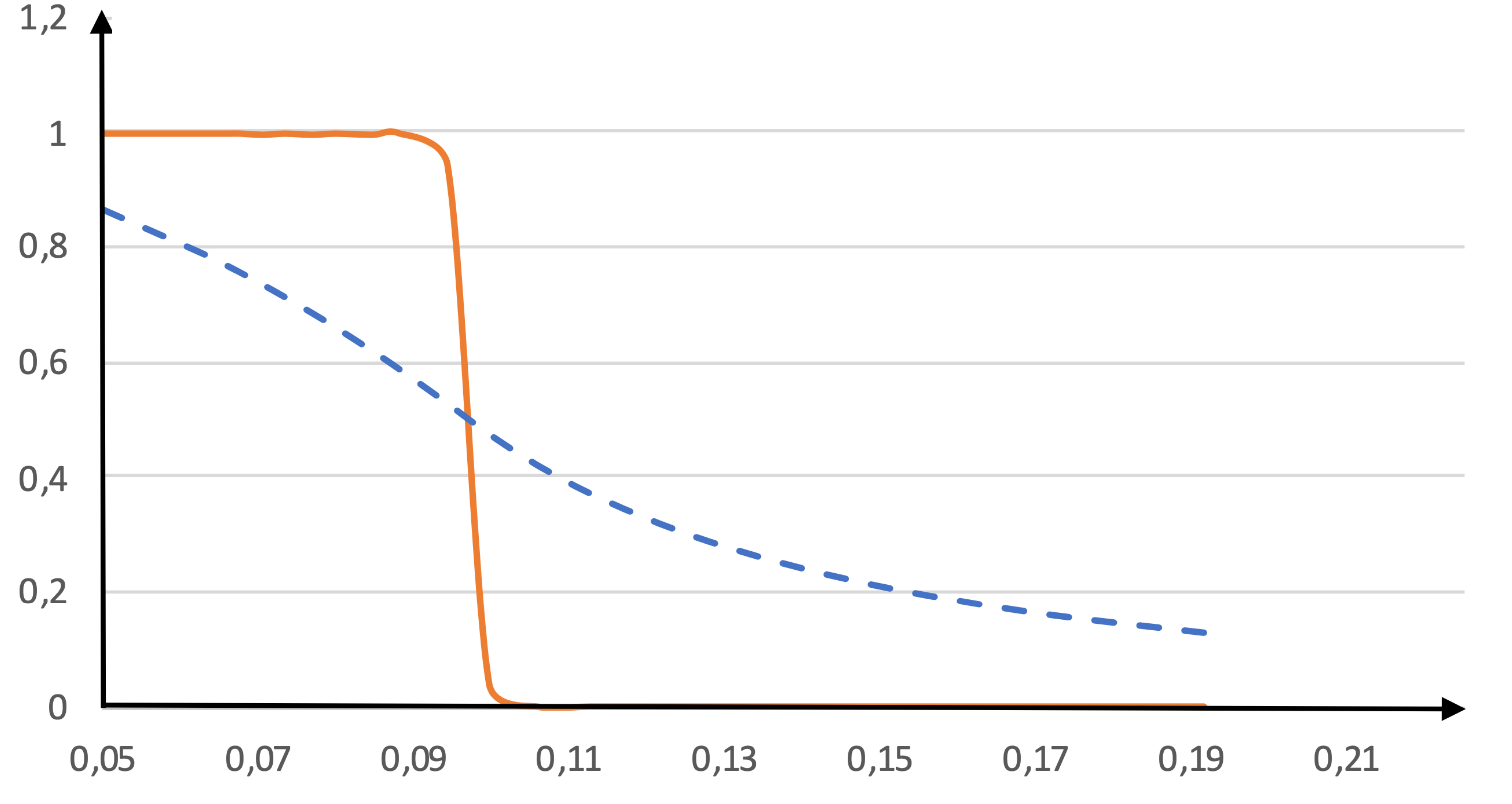}
}
&
\hspace*{-0.5cm}
{
\includegraphics[width=9.cm]{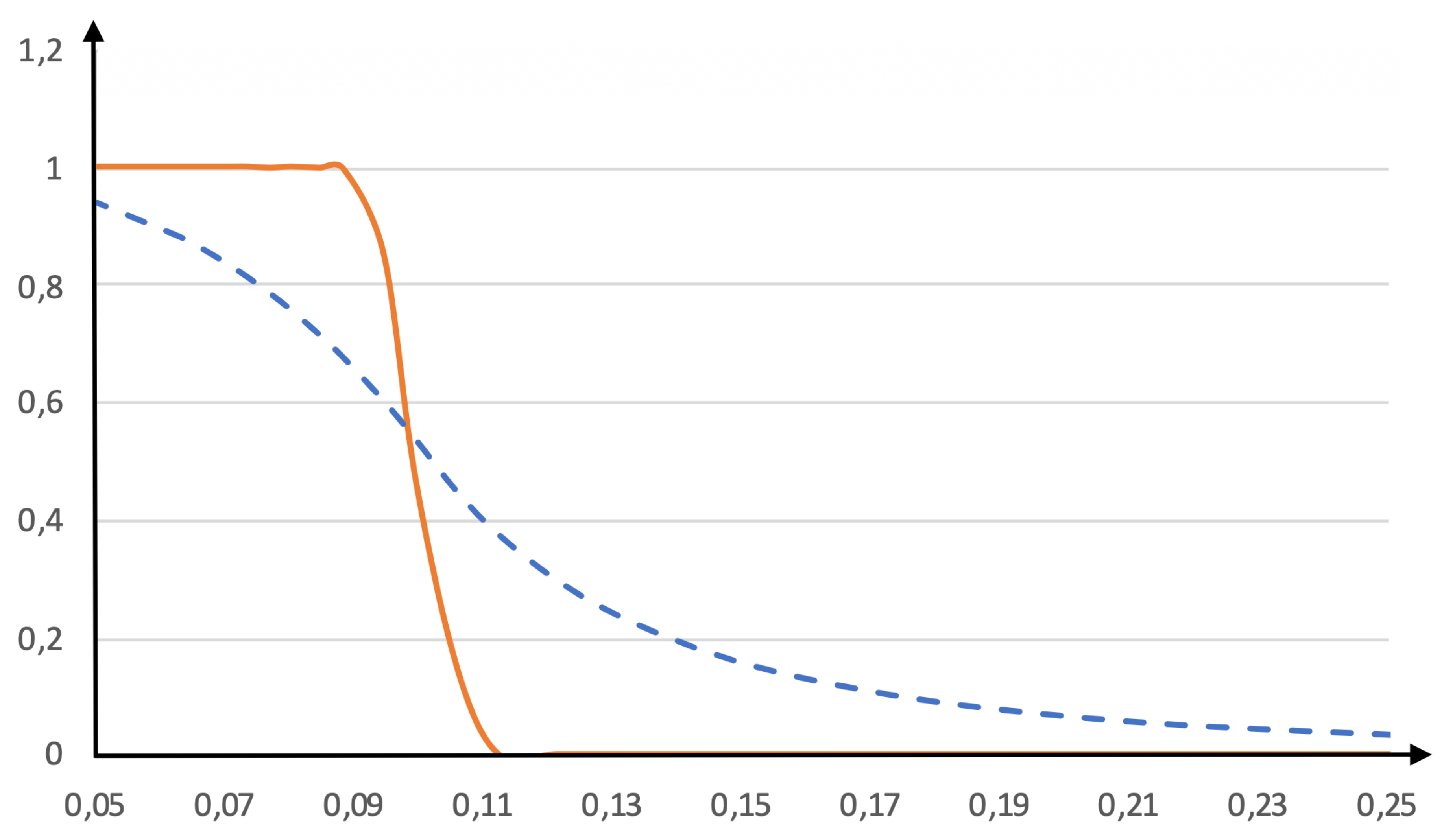}
}
\end{tabular}
\caption{$P_2$ versus $P_4$ (left) and $P_1$ versus $P_4$ (right) for the Runge function with $\alpha=2000$. Comparison between the statistical frequencies (full) and the probabilistic law (\ref{Nonlinear_Prob}) (dotted). 500 meshes are used for each value of $h$.} \label{P1-P2-P4-sigmoid}
\end{figure}

\begin{figure}[htb]
  \centering
  \includegraphics[width=9.cm]{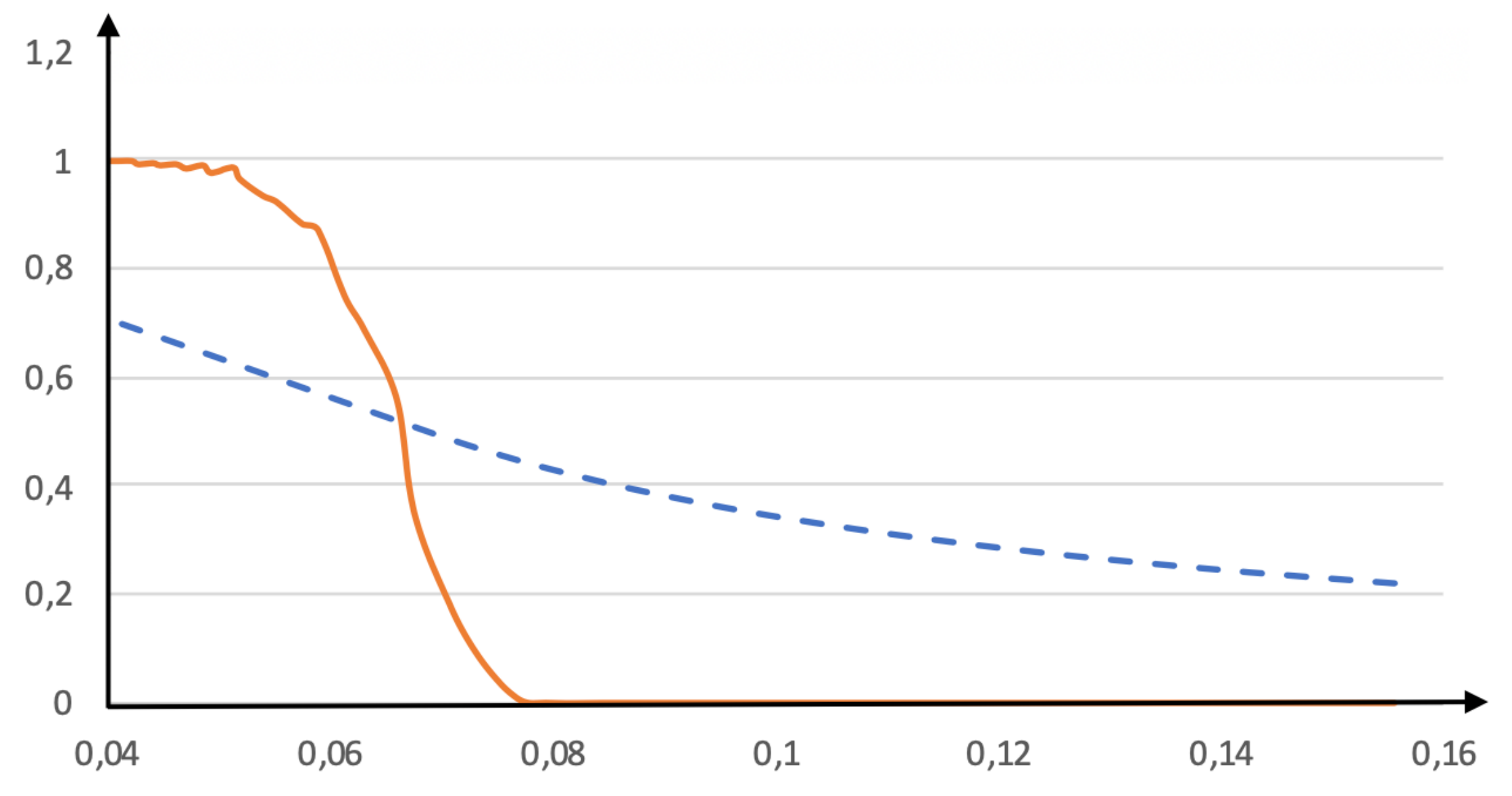}
  \caption{$P_2$ versus $P_3$ for the  Runge function with $\alpha=2000$. Comparison between the statistical frequencies (full) and the probabilistic law (\ref{Nonlinear_Prob}) (dotted). 500 meshes are used for each value of $h$.}
 \label{P2-P3-Runge-sigmoid}
\end{figure}

\noindent Next, we numerically assessed the validity of the present approach when finite element $P_4$ is involved. We first compared $P_{3}$ with $P_{4}$, then $P_{2}$ with $P_{4}$. To this end, we used the Runge function with $\alpha=2000$, and we also considered $N=500$ different meshes. The computed value of $\widehat{\,h^{*}_{k,m}}$ we obtained are $\widehat{\,h^{*}_{3,4}}\simeq 0.24$ for $P_3$-$P_4$ and $\widehat{\,h^{*}_{2,4}}\simeq 0.094$ for $P_2$-$P_4$. The results are depicted in Fig. \ref{P3P4-P2P4-alpha-2000} and show, like previously, that the statistical frequencies behave very similarly to the two-steps probabilistic law (\ref{Heaviside_Prob}).\sa
\noindent The last illustration of this subsection is devoted to the comparison between the statistical frequencies and the sigmoid probabilistic law defined in (\ref{Nonlinear_Prob}). We considered comparisons between finite elements $P_2$ and $P_4$, then between $P_1$ and $P_4$. \sa
We followed the same procedure as above, again with the same parameters ($N=500$ and $\alpha=2000$). The results are depicted in Fig. \ref{P1-P2-P4-sigmoid}. As one can see, there is a weaker fitting between the two curves than with the two-steps law (\ref{Heaviside_Prob}), even if the trend is still correct. Remark also that the fit is better in the $P_1$-$P_4$ case than in the $P_2$-$P_4$ case. More generally, the greater the $m-k$ difference, the better the match. Hence, the sigmoid model also gave a correct trend, but was less precise and satisfying than the two-steps law, in particular when $m-k=1$, for instance when one compares $P_{2}$ with $P_{3}$, see Fig. \ref{P2-P3-Runge-sigmoid}. Indeed, in that case, the first part of (\ref{Nonlinear_Prob}) is a linear decreasing function of $h$ (for any given fixed $h^{*}_{k,m}$), and the second one decreases like $1/h$. However, if difference $m-k=2$, for instance when comparing $P_{2}$ to $P_{4}$, the fit is better. Indeed, the first part of (\ref{Nonlinear_Prob}) is a decreasing function $\simeq -h^{2}$ (for any given fixed $h^{*}_{k,m}$), and the second one decreases like $1/h^{2}$. So, depending on the difference $m-k$, the sigmoid law remains to some extent relevant, where high order finite element (with $m$ around $20$-$25$) are sometimes used \cite{Mitc15}.\sa
This shows that the two-steps model works well, but is a bit ``rough'' (essentially binary), whereas the sigmoid law is probably too ``rigid'' and has to be make more ``flexible'' to obtain a better fit with the statistical results. For this reason, we are currently working on a more general approach which corresponds to a new generation of probabilistic laws that better fit the statistical frequencies. \sa

\subsection{A smooth example}
\begin{figure}[htbp!]
\begin{tabular}{lr}
{
\hspace*{-0.5cm}
\includegraphics[width=8.cm]{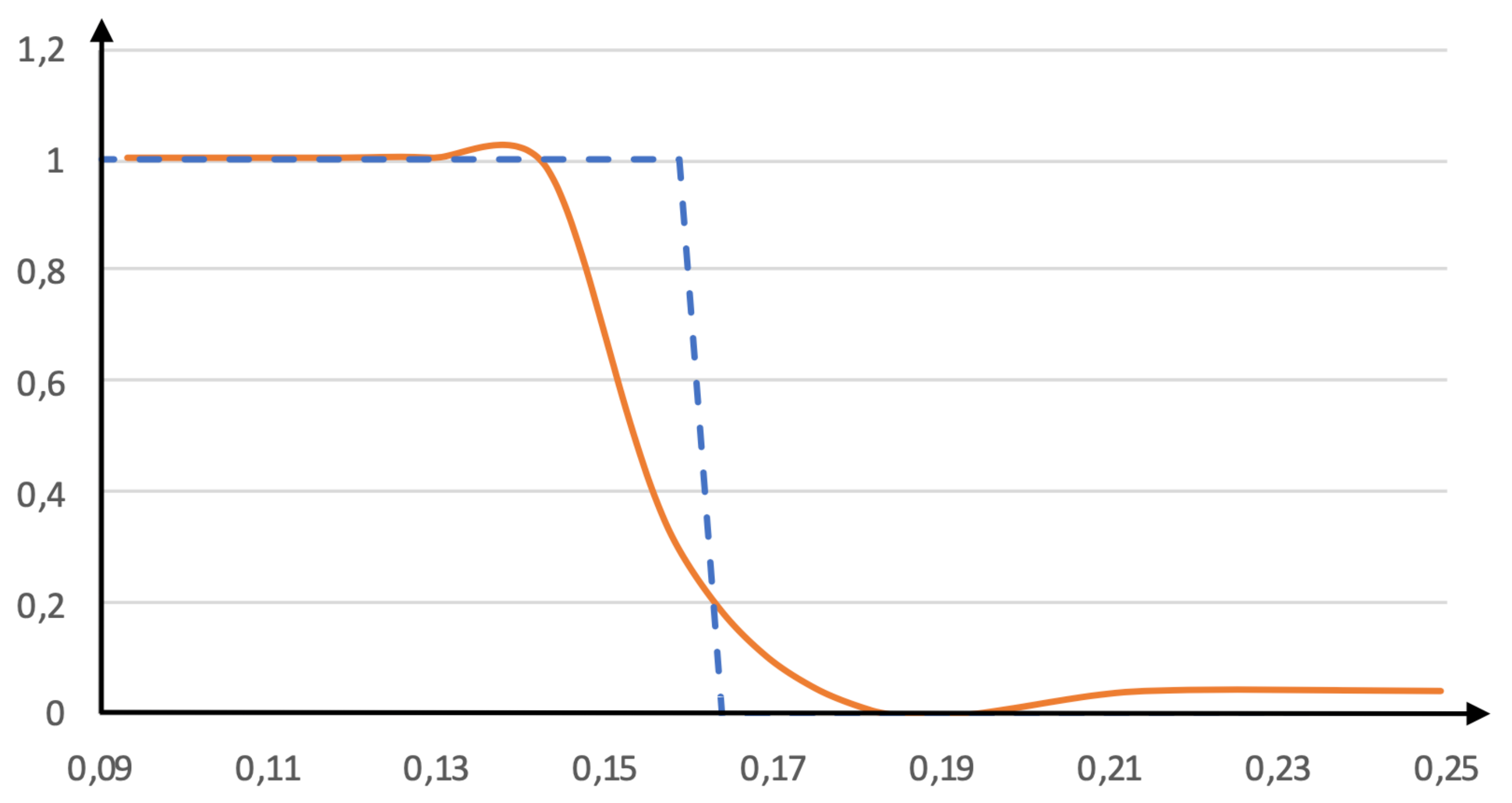}
}
&
\hspace*{-0.5cm}
{
\includegraphics[width=8.cm]{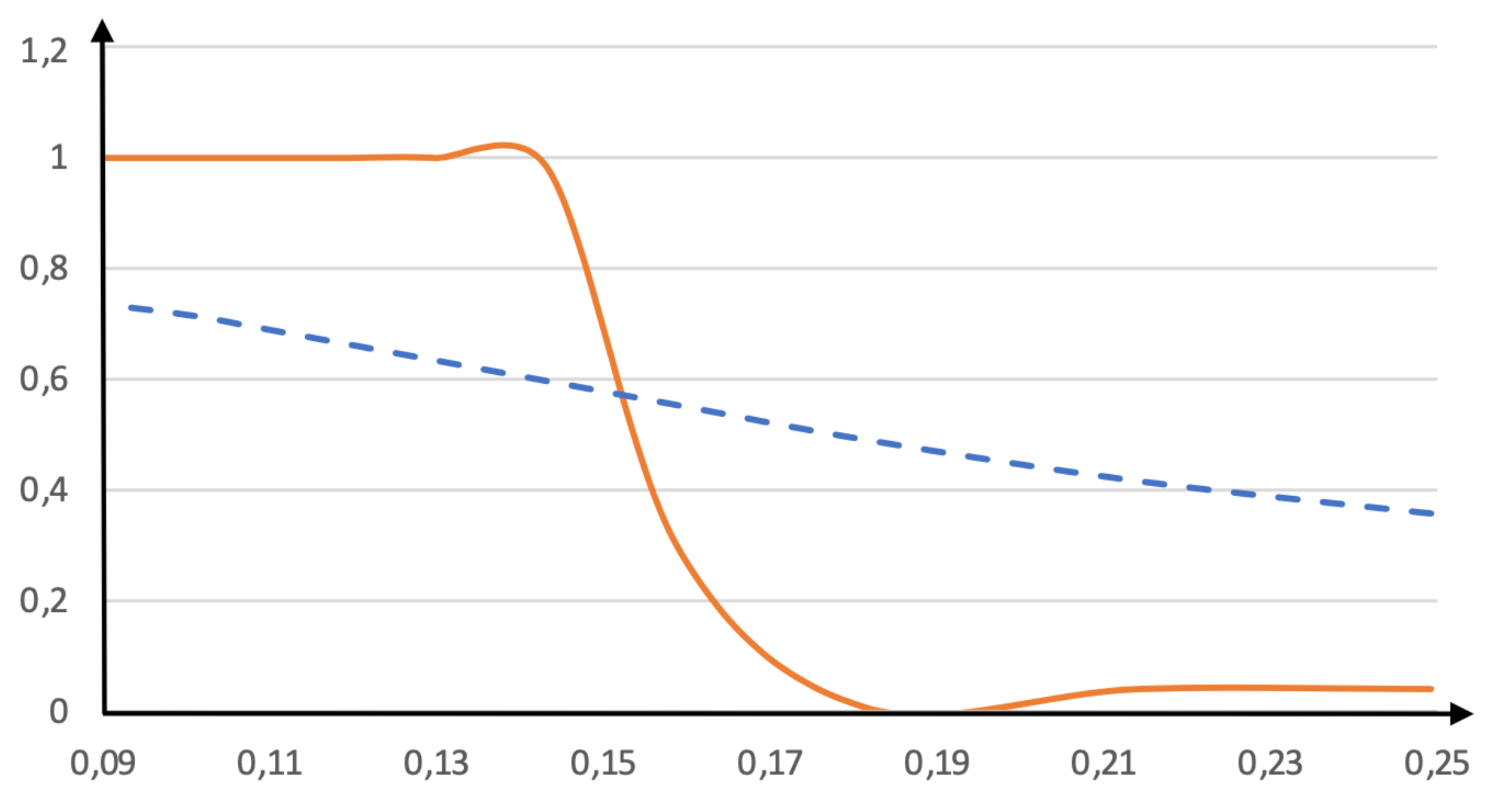}
}
\end{tabular}
\caption{$P_2$ versus $P_3$ for the smooth case. Comparison between the statistical frequencies (full) and the two probabilistic laws (dotted). 500 meshes are used for each value of $h$. Left: comparison with the two-steps probabilistic law (\ref{Heaviside_Prob}) - Right: comparison with the sigmoid probabilistic law (\ref{Nonlinear_Prob}).} \label{P2-P3-regular}
\end{figure}

\noindent In this subsection, we illustrate the probabilistic laws for a very smooth solution to the variational  problem (\ref{VP}). To build such a case, we chose $q=2\pi^2 \sin(\pi x) \cos(\pi y)$ in (\ref{edpnum}), so that $u(x,y)=\sin(\pi x) \cos(\pi y)$ is the exact solution of the problem, provided that the Dirichlet boundary condition $h$ is taken as the trace of $u(x,y)$ on the boundary $\partial \Omega$. This can be written as:
$$
\left\{
\begin{array}{ll}
h(x,0) = \sin(\pi x)   & h(0,y) = 0, \\
h(x,1) = -\sin(\pi x),  & h(1,y) = 0.
\end{array}
\right.
$$
As previously, we first compute $\widehat{\,h^{*}_{k,m}}$ defined by (\ref{hstar_estime}), then we compute the probabilistic models introduced above. After that, we compare these results to the statistical frequencies.\sa
For example, we consider the finite elements $P_2$ and $P_3$: we depicted in Fig. \ref{P2-P3-regular} (left) the statistical frequencies and the probabilistic law (\ref{Heaviside_Prob}). As in all the other numerical experiments, 500 meshes have been used for each value of $h$, where we found a value of $\widehat{\,h^{*}_{2,3}}$ approximately equal to $0.18$. As one can see, even in this case, there is a good agreement between the statistical frequencies and the probabilistic law (\ref{Heaviside_Prob}). However, the comparison with the sigmoid law (\ref{Nonlinear_Prob}) (right part of Fig. \ref{P2-P3-regular}) gave only a global trend and was not really accurate. Here again, as for the Runge example, it will be improved by the above-mentioned new generation of probabilistic laws.\sa
\section{Conclusion}
\noindent In this paper, we proposed to apply the probabilistic approach we developed in \cite{CMAM_2019} to numerical examples. It enabled us to evaluate the relative accuracy between two Lagrange finite elements $P_k$ and $P_m, (k<m)$, for a fixed value of the mesh size $h$. Our approach, which is based on a geometrical interpretation of the error estimate, considers the approximation errors as random variables. Two probabilistic laws were derived, a so-called "two steps" law and  a "sigmoid" one, depending on the probabilistic assumptions which were made on the corresponding random variables.\sa
For the finite elements we considered, we illustrated, using several examples, the property that, depending on the position of $h$ with respect to the critical value $h^{*}_{k,m}$, we can actually estimate which of finite elements $P_k$ and $P_m$ is more likely accurate. This overturns the common misconception that finite elements $P_m$ are always more precise than $P_k$ if $m > k$, regardless of the mesh size $h$. In particular, this shows cases where a $P_m$ finite element \emph{surely} is overqualified. As a consequence, a significant reduction of implementation time and execution cost can be obtained without loss of accuracy. Such a phenomenon was already observed by using data-mining techniques (see \cite{AsCh11}, \cite{AsCh13}, \cite{AsCh16} and \cite{AsCh17}).\sa
However if, in the proposed examples, the first investigated law (the two-steps law) fit the numerical results satisfactorily, the second proposed law (the sigmoid one) produces only a trend and is not accurate enough. Indeed, the results show that the statistical  frequencies behave similarly to the two-steps probabilistic law, in both the smooth and stiff examples. However, there is a weaker fitting between the statistical error and the sigmoid law, particularly when difference $m-k$ is small. To address this issue, we are currently working on a new probabilistic framework which corrects the gap between the statistics and a ``generalized''  probability law.\sa
Finally, note that this approach is not limited to finite element methods, and can be generalized to other approximation methods: given several different numerical methods and their error estimates, it would be possible to order them by evaluating which is the most probably accurate.\sa
\textbf{\underline{Homages}:} The author wants to warmly dedicate this research to pay homage to the memory of Professors Andr\'e Avez and G\'erard Tronel, who broadly promoted the passion of research and teaching in Mathematics.
\end{document}